\documentstyle[12pt,thmsa,sw20lart]{article}

\input{tcilatex}
\begin{document}

\author{TEODOR OPREA \\
University of Bucharest,\\
Faculty of Mathematics and Informatics\\
$14$ Academiei St., code 010014, Bucharest, Romania\\
e-mail: teodoroprea@yahoo.com}
\title{ON A GEOMETRIC\ INEQUALITY }
\date{}
\maketitle

\bigskip {\bf Abstract.} As we showed in [3], a geometric inequality can be
regarded as an optimization problem. In this paper we find another proof for
a Chen's inequality, regarding the Ricci curvature [2] and we improve this
inequality in the Lagrangian case.

\medskip

2000{\bf \ }{\it Mathematics subject classification: }53C21, 53C24, 53C25,
49K35.


\bigskip

{\bf 1. Optimizations on Riemannian submanifolds. }Let $(N,\widetilde{g})$
be a Riemannian manifold, $(M,g)$ a Riemannian submanifold of it, and $f\in 
{\cal F}(N).$ To these ingredients we attach the optimum problem\newline

\[
(1)\text{ }\min\limits_{x\in M}f(x). 
\]

Let's remember the result obtained in [3].

\ 

T{\footnotesize HEOREM} 1.1.{\sl \ If }$x_{0}\in M${\sl \ is the solution of
the problem }$(1)${\sl , then}

{\sl \ }

i){\sl \ }$($grad f$)(x_{0})${\sl \ }$\in T_{x_{0}}^{\perp }M,$

{\sl \ }

ii){\sl \ the bilinear form }

\[
\begin{array}{c}
\alpha {\sl \ :\ }T_{x_{0}}M\times T_{x_{0}}M{\sl \ }\rightarrow {\sl \ }R,
\\ 
\alpha (X,Y)=\text{Hess}_{f}(X,Y)+\widetilde{g\text{ }}(h(X,Y),(\text{grad f}%
)(x_{0}))
\end{array}
\]
{\sl is positive semidefinite, where }$h${\sl \ is the second fundamental
form of the submanifold }$M${\sl \ in }$N.$

\medskip

R{\footnotesize EMARK}{\bf . }The bilinear form $\alpha $ is nothing else
but Hess$_{\left. f\right| M}$.

\medskip

A very nice application of this result is the next problem: {\sl find the
2-plane in the tangent space at the given point }$x$ {\sl of a Riemannian
manifold (}$M,g)$ {\sl which minimize the sectional curvature.} A equivalent
conditioned extremum problem is

\[
(2)\min \text{ }R(X,Y,X,Y), 
\]
\[
\text{subject to }\parallel X\parallel =1,\parallel Y\parallel =1\text{, }%
g(X,Y)=0, 
\]
where $X$ and $Y$ are two vectors from $T_{x}M$.

\medskip

{\sl A 2-plane }$\pi \subset T_{x}M${\sl , }$\pi =$Sp$\{X,Y\}${\sl \ which
verifies the first condition from Theorem 1.1 is called critical plane for
the sectional curvature at the point }$x$. Using Theorem 1.1, in [3] we
showed that {\sl a 2-plane }$\pi ${\sl \ is a critical plane for the
sectional curvature at the point }$x${\sl \ if and only if for every tangent
vectors }$U$, $V$, $W\in $ $\pi $ {\sl the vector }$R(U,V)W${\sl \ lies in} $%
\pi $, {\sl where }${\sl R}$ {\sl is} {\sl the curvature tensor} {\sl of the
Riemannian manifold (}$M,g${\sl )}.

\bigskip

{\bf 2. The Ricci curvature of a submanifold in a real space form.} In this
section we give another demonstration for the next inequality:

\medskip

T{\footnotesize HEOREM} 2.1(Chen [2]). {\sl Let }$M${\sl \ be a }$n-${\sl %
dimensional Riemannian submanifold of a real space form (}$\widetilde{M}%
(c),g)${\sl \ and }$x${\sl \ a point in }$M.${\sl \ Then, for each unit
vector }$X\in T_{x}M,${\sl \ we have} 
\[
\text{Ric}(X)\leq (n-1)c+\frac{n^{2}}{4}\left\| H\right\| ^{2}, 
\]
{\sl where }$H${\sl \ is the mean curvature vector of }$M$ {\sl in} $%
\widetilde{M}(c)${\sl \ and }Ric{\sl (}$X${\sl ) the Ricci curvature of }$M$%
{\sl \ at }$x${\sl .}

\medskip

{\it Proof.} We fix the point $x$ in $M,$the vector $X\in T_{x}M$, with%
\newline
$\left\| X\right\| =1$, the orthonormal frame $\{e_{1},e_{2},...,e_{n}\}$ in 
$T_{x}M$ such that $e_{1}=X$ and $\{e_{n+1},e_{n+2},...,e_{m}\}$ a
orthonormal frame in $T_{x}^{\perp }M$ .

\medskip

From Gauss equation we have\newline
(1) $\widetilde{R}(e_{1},e_{j},e_{1},e_{j})=R(e_{1},e_{j},e_{1},e_{j})-%
\widetilde{g}(h(e_{1},e_{1}),h(e_{j},e_{j}))+$\newline
$+\widetilde{g}(h(e_{1},e_{j}),h(e_{1},e_{j}))=R(e_{1},e_{j},e_{1},e_{j})-%
\dsum\limits_{r=n+1}^{m}(h_{11}^{r}h_{jj}^{r}-(h_{1j}^{r})^{2}),$ $j\in 
\overline{2,n}.$

Using the fact that the sectional curvature of $\widetilde{M}(c)$ is
constant, we obtain\newline
(2) $(n-1)c=$Ric$(X)-\dsum\limits_{r=n+1}^{m}\dsum%
\limits_{j=2}^{n}(h_{11}^{r}h_{jj}^{r}-(h_{1j}^{r})^{2})$, therefore\newline
(3) Ric$(X)-(n-1)c=\dsum\limits_{r=n+1}^{m}\dsum%
\limits_{j=2}^{n}(h_{11}^{r}h_{jj}^{r}-(h_{1j}^{r})^{2})\leq
\dsum\limits_{r=n+1}^{m}\dsum\limits_{j=2}^{n}h_{11}^{r}h_{jj}^{r}.$

\ 

For $r\in \overline{n+1,m}$, let us consider the quadratic form

\[
\begin{array}{c}
f_{r}:R^{n}\rightarrow R, \\ 
f_{r}(h_{11}^{r},h_{22}^{r},...,h_{nn}^{r})=\dsum%
\limits_{j=2}^{n}h_{11}^{r}h_{jj}^{r}
\end{array}
\]
and the constrained extremum problem 
\[
\max f_{r}\text{,} 
\]

\[
\text{subject to }P:h_{11}^{r}+h_{22}^{r}+...+h_{nn}^{r}=k^{r}, 
\]
where $k^{r}$ is a real constant.

\bigskip \medskip

The gradient vector of $f_{r}$ is given by

\[
\text{grad}f_{r}=(\dsum%
\limits_{j=2}^{n}h_{jj}^{r},h_{11}^{r},...,h_{11}^{r}). 
\]

\ 

Let us denote with $p=(h_{11}^{r},h_{22}^{r},...,h_{nn}^{r})$ a solution of
the problem in question.

\ 

As grad$f_{r}$ is normal to $P$ at the point $p$, we obtain \newline
(4) $h_{11}^{r}=\dsum\limits_{j=2}^{n}h_{jj}^{r}=\frac{k^{r}}{2}$ .

Let $q\in P$ be an arbitrary point.

The bilinear form $\alpha :T_{q}P\times T_{q}P\rightarrow R$ has the
expression 
\[
\alpha (X,Y)=\text{Hess}(f_{r})(X,Y)+\newline
\langle h^{^{\prime }}(X,Y),\text{grad}f_{r}(q)\rangle , 
\]
where $h^{^{\prime }}$ is the second fundamental form of $P$ in $R^{n}$, and 
$\langle $ , $\rangle $ is the standard inner-product on $R^{n}$.

\medskip

In the standard frame of $R^{n}$, the Hessian of $f_{r}$ has the matrix

\ 

\[
\left( 
\begin{array}{cccc}
0 & 1 & . & 1 \\ 
1 & 0 & . & 0 \\ 
. & . & . & . \\ 
1 & 0 & . & 0
\end{array}
\right) 
\]
$.$

\ 

A vector $X\in T_{q}P$ satisfies $\sum\limits_{i=1}^{n}X^{i}=0$.

\ 

As $P$ is totally geodesic in $R^{n},$ we have $\alpha
(X,X)=2\sum\limits_{j=2}^{n}X^{1}X^{j}=$\newline
$=(X^{1}+\sum\limits_{j=2}^{n}X^{i})^{2}-(\sum%
\limits_{j=2}^{n}X^{j})^{2}-(X^{1})^{2}=-(\sum%
\limits_{j=2}^{n}X^{j})^{2}-(X^{1})^{2}\leq 0.$ So $\left. f_{r}\right| P$
is a convex function, therefore the points which satisfies the relation (4)
are global maximum points for $\left. f_{r}\right| P$.

\ 

One gets\newline
(5) $f_{r}\leq \frac{(k^{r})^{2}}{4}=$ $\frac{1}{4}(\dsum%
\limits_{i=1}h_{ii}^{r})^{2}=\frac{n^{2}}{4}(H^{r})^{2}.$

\ 

Using (3) and (5) we find\newline
(6) Ric$(X)-(n-1)c\leq \dsum\limits_{r=n+1}^{m}\frac{n^{2}}{4}(H^{r})^{2}=%
\frac{n^{2}}{4}\left\| H\right\| ^{2}$.

\bigskip

{\bf 3. The Ricci curvature of a Lagrangian submanifold in a complex space
form.} Let $(\widetilde{M},\widetilde{g},J)$ be a K\"{a}hler manifold of
real dimension $2m.$ A submanifold $M$ of dimension $n$ of $(\widetilde{M},%
\widetilde{g},J)$ is called a totally real submanifold if for any point $x$
in $M$ the relation $J(T_{x}M)\subset T_{x}^{\perp }M$ holds.

If, in addition, $n=m,$ then $M$ is called Lagrangian submanifold. For a
Lagrangian submanifold, the relation $J(T_{x}M)=T_{x}^{\perp }M$ occurs.

A K\"{a}hler manifold with constant holomorphic sectional curvature is
called a complex space form and is denoted by $\widetilde{M}(c)$. The
Riemann curvature tensor $\widetilde{R}$ of $\widetilde{M}(c)$ satisfies the
relation

$\widetilde{R}(X,Y)Z=\frac{c}{4}\{\widetilde{g}(Y,Z)X-\widetilde{g}(X,Z)Y+%
\widetilde{g}(JY,Z)JX-\widetilde{g}(JX,Z)JY+$\newline
$+2\widetilde{g}(X,JY)JZ\}.$

A totally real submanifold of real dimension $n$ in a complex space form $%
\widetilde{M}(c)$ of real dimension $2m$ verifies a Chen's inequality:

\medskip

T{\footnotesize HEOREM} 3.1(Chen). {\sl Let }$M${\sl \ be a }$n-${\sl %
dimensional Riemannian submanifold of a complex space form (}$\widetilde{M}%
(c),g)${\sl \ and }$x${\sl \ a point in }$M.${\sl \ Then, for each unit
vector }$X\in T_{x}M,${\sl \ we have} 
\[
\text{Ric}(X)\leq (n-1)\frac{c}{4}+\frac{n^{2}}{4}\left\| H\right\| ^{2}, 
\]
{\sl where }$H${\sl \ is the mean curvature vector of }$M$ {\sl in} $%
\widetilde{M}(c)${\sl \ and }Ric{\sl (}$X${\sl ) the Ricci curvature of }$M$%
{\sl \ at }$x${\sl .}

\medskip

R{\footnotesize EMARK}{\bf . }i) If $M$ is a totally real submanifold of
real dimension $n$ in a complex space form $\widetilde{M}(c)$ of real
dimension $2m,$ then

\begin{center}
\medskip $A_{JY}X=-Jh(X,Y)=A_{JX}Y,$ $\forall $ $X,$ $Y\in {\cal X}(M).$
\end{center}

ii) Let $m=n$ ($M$ is Lagrangian in $\widetilde{M}(c)$). If we consider the
point $x\in M$, the orthonormal frames $\{e_{1},...,e_{n}\}$ in $T_{x}M$ and 
$\{Je_{1},...,Je_{n}\}$ in $T_{x}^{\perp }M$, then

\medskip 
\[
h_{jk}^{i}=h_{ik}^{j},\forall \text{ }i,j,k\in \overline{1,n}, 
\]
where $h_{jk}^{i}$ is the component after $Je_{i}$ of the vector $%
h(e_{j},e_{k}).$

\ 

With these ingredients we prove the next result which improve Chen's
inequality in the Lagrangian case.

\medskip

T{\footnotesize HEOREM} 3.2.{\bf \ }{\sl Let }$M${\sl \ be a Lagrangian
submanifold in a complex space form }$\widetilde{M}(c)${\sl \ of real
dimension }$2n${\sl , }$n\geq 2$, $x$ {\sl a point in} $M$ {\sl and} $X$ 
{\sl a unit tangent vector in }$T_{x}M.${\sl \ Then we have}

\[
\text{Ric}(X)\leq \frac{n-1}{4}(c+n\left\| H\right\| ^{2}). 
\]

{\it Proof}. We fix the point $x$ in $M,$ the tangent vector $X\in T_{x}M$,
with $\left\| X\right\| =1$, the orthonormal frame $\{e_{1},e_{2},...,e_{n}%
\} $ in $T_{x}M$ such that $e_{1}=X$ and $\{Je_{1},Je_{2},...,Je_{n}\}$ a
orthonormal frame in $T_{x}^{\perp }M$ .

\ 

From the Gauss equation we get\newline
(1) $\widetilde{R}(e_{1},e_{j},e_{1},e_{j})=R(e_{1},e_{j},e_{1},e_{j})-%
\widetilde{g}(h(e_{1},e_{1}),h(e_{j},e_{j}))+$\newline
$+\widetilde{g}(h(e_{1},e_{j}),h(e_{1},e_{j}))=R(e_{1},e_{j},e_{1},e_{j})-%
\dsum\limits_{r=1}^{n}(h_{11}^{r}h_{jj}^{r}-(h_{1j}^{r})^{2})$, $\forall $ $%
j\in \overline{2,n}.$ Therefore\newline
(2) $(n-1)\frac{\text{{\normalsize c}}}{4}=$Ric$(X)-\dsum\limits_{r=1}^{n}%
\dsum\limits_{j=2}^{n}(h_{11}^{r}h_{jj}^{r}-(h_{1j}^{r})^{2})$, which implies%
\newline
(3) Ric$(X)-(n-1)\frac{\text{{\normalsize c}}}{4}=\dsum\limits_{r=1}^{n}%
\dsum\limits_{j=2}^{n}(h_{11}^{r}h_{jj}^{r}-(h_{1j}^{r})^{2})\leq
(\dsum\limits_{r=1}^{n}\dsum\limits_{j=2}^{n}h_{11}^{r}h_{jj}^{r})-\dsum%
\limits_{j=2}^{n}(h_{1j}^{1})^{2}-$\newline
$-\dsum\limits_{j=2}^{n}(h_{1j}^{j})^{2}=(\dsum\limits_{r=1}^{n}\dsum%
\limits_{j=2}^{n}h_{11}^{r}h_{jj}^{r})-\dsum%
\limits_{j=2}^{n}(h_{11}^{j})^{2}-\dsum\limits_{j=2}^{n}(h_{jj}^{1})^{2}.$

\ 

Let us consider the quadratic forms 
\[
\begin{array}{c}
f_{1},\text{ }f_{r}:R^{n}\rightarrow R\text{, } \\ 
f_{1}(h_{11}^{1},h_{22}^{1},...,h_{nn}^{1})=\dsum%
\limits_{j=2}^{n}h_{11}^{1}h_{jj}^{1}-\dsum\limits_{j=2}^{n}(h_{jj}^{1})^{2},
\\ 
f_{r}(h_{11}^{r},h_{22}^{r},...,h_{nn}^{r})=\dsum%
\limits_{j=2}^{n}h_{11}^{r}h_{jj}^{r}-(h_{11}^{r})^{2},
\end{array}
\]
where $r\in \overline{2,n}.$

We need the maximum of $f_{1}$ and $f_{2}$. For $f_{r}$, $r\in \overline{3,n}
$, we can solve similar problems.

We start with the problem

\[
\max f_{1}, 
\]
\[
\text{cu restric\c{t}ia }P:h_{11}^{1}+h_{22}^{1}+...+h_{nn}^{1}=k^{1}, 
\]
where $k^{1}$ is a real constant.

\ 

The partial derivatives of the function $f_{1}$ are\newline
(4) $\frac{\partial f_{1}}{\partial h_{11}^{1}}=\dsum%
\limits_{j=2}^{n}h_{jj}^{1}$,\newline
(5) $\frac{\partial f_{1}}{\partial h_{jj}^{1}}=h_{11}^{1}-2h_{jj}^{1}$, $%
\forall $ $j\in \overline{2,n}$.

\ 

As for a solution $(h_{11}^{1},h_{22}^{1},...,h_{nn}^{1})$ of the problem in
question the vector grad$f_{1}$ is normal to $P$, we obtain\newline
(6) $h_{22}^{1}=h_{33}^{1}=...=h_{nn}^{1}=a^{1}$ and\newline
(7) $h_{11}^{1}-2h_{22}^{1}=\dsum\limits_{j=2}^{n}h_{jj}^{1}.$

\ 

Using (6) and (7) we find\newline
(8) $h_{11}^{1}=(n+1)a^{1}.$

\ 

From the relation $h_{11}^{1}+h_{22}^{1}+...+h_{nn}^{1}=k^{1},$ we get%
\newline
(9) $(n+1)a^{1}+(n-1)a^{1}=k^{1}$, therefore\newline
(10) $a^{1}=\frac{k^{1}}{2n}$.

\medskip

As $f_{1\text{ }}$is obtained from the function studied in previous section
by subtracting some square terms, $f_{1}\left| P\right. $ will have the
Hessian negative definite. Consequently the point $%
(h_{11}^{1},h_{22}^{1},...,h_{nn}^{1})$ given by the relations (6), (8) and
(10) is a maximum point, and hence\newline
(11) $f_{1}\leq (n+1)a^{1}(n-1)a^{1}-(n-1)(a^{1})^{2}=$\newline
$=n(n-1)(a^{1})^{2}=\frac{n-1}{4n}(k^{1})^{2}=\frac{n(n-1)}{4}(H^{1})^{2}.$

\medskip

Further on, we shall consider the problem

\[
\max f_{2}, 
\]
\[
\text{subject to }P:h_{11}^{2}+h_{22}^{2}+...+h_{nn}^{2}=k^{2}, 
\]
where $k^{2}$ is a real constant.

\ 

\ 

The first two partial derivatives of the function $f_{2}$ are\newline
(13) $\frac{\partial f_{2}}{\partial h_{11}^{2}}=\dsum%
\limits_{j=2}^{n}h_{jj}^{2}-2h_{11}^{2}$ ,\newline
(14) $\frac{\partial f_{2}}{\partial h_{jj}^{2}}=h_{11}^{2},$ $\forall $ $%
j\in \overline{2,n}$.

\ 

As for a solution $(h_{11}^{2},h_{22}^{2},...,h_{nn}^{2})$ of the problem in
question the vector grad$f_{2}$ is normal to $P$, we obtain\newline
(15) $3h_{11}^{2}=\dsum\limits_{j=2}^{n}h_{jj}^{2}=3a^{2}$.

\ 

Using the relation $h_{11}^{2}+h_{22}^{2}+...+h_{nn}^{2}=k^{2},$ we obtain%
\newline
(16) $a^{2}+3a^{2}=k^{2}$, therefore\newline
(17) $a^{2}=\frac{k^{2}}{4}.$

\ 

With an similar argument to those in the previous problem we obtain that the
point $(h_{11}^{2},h_{22}^{2},...,h_{nn}^{2})$ given by the relations (15)
and (17) is a maximum point. Therefore\newline
\newline
(18) $f_{2\text{ }}\leq a^{2}3a^{2}-(a^{2})^{2}=2(a^{2})^{2}=\frac{%
(k^{2})^{2}}{8}=\frac{n^{2}}{8}(H^{2})^{2}.$

\ 

Similarly one gets\newline
(19) $f_{r}\leq \frac{n^{2}}{8}(H^{r})^{2}$, $\forall $ $r\in \overline{2,n}%
. $

\ 

As $\frac{n(n-1)}{4}\geq \frac{n^{2}}{8}$, $\forall $ $n\geq 2$, using (11)
and (19) we find\newline
(20) $f_{r}\leq \frac{n(n-1)}{4}(H^{r})^{2}$, $\forall $ $r\in \overline{1,n}%
.$

\ 

From (3) and (20) it follows\newline
(21) Ric$(X)-(n-1)\frac{c}{4}\leq \frac{n(n-1)}{4}\dsum%
\limits_{r=1}^{n}(H^{r})^{2}=\frac{n(n-1)}{4}\left\| H\right\| ^{2}$,
therefore\newline
(22) Ric$(X)\leq \frac{n-1}{4}(c+n\left\| H\right\| ^{2}).$

\end{document}